\documentclass[letterpaper, 10pt, conference]{ieeeconf}
\IEEEoverridecommandlockouts 
\usepackage{amssymb}
\usepackage{amsmath}

\usepackage{amsthm}
\usepackage{cite}
\usepackage{graphics}
\usepackage{hyperref}
\usepackage{multirow}
\usepackage{tikz}

\newtheoremstyle{named}{}{}{}{}{\bfseries}{.}{.5em}{\thmnote{#3} Principle}
\theoremstyle{named}

\theoremstyle{definition}
\newtheorem{definition}{Definition}

\def\u{\mathbf{u}}
\def\w{\mathbf{w}}
\def\x{\mathbf{x}}
\def\K{\mathbf{K}}

\def\R{\mathbf{\Phi}_x}
\def\M{\mathbf{\Phi}_u}

\newif\ifshowWriterComment
\newcommand\writercomment[3]{\expandafter\newcommand\csname #2\endcsname[1]{\ifshowWriterComment{\color{#3} #1: ##1}\fi}}

\writercomment{Lisa}{Lisa}{purple}

\showWriterCommenttrue

\title{\LARGE \bf
	Descending Predictive Feedback: \\ From Optimal Control to the Sensorimotor System
}

\author{Jing Shuang (Lisa) Li, Anish A. Sarma, and John C. Doyle
	\thanks{J. S. Li and J. C. Doyle are with Computing and Mathematical Sciences, California Institute of Technology. A. A. Sarma is with Computation and Neural Systems, California Institute of Technology. {\tt\small \{jsli, aasarma, doyle\}@caltech.edu},
	}%
}
{\tiny }
\begin{document}	
	\maketitle
	\thispagestyle{empty}
	\pagestyle{empty}
	
	\begin{abstract}

Descending predictive feedback (DPF) is an ubiquitous yet unexplained phenomenon in the central nervous system. Motivated by recent observations on motor-related signals in the visual system, we approach this problem from a sensorimotor standpoint and make use of optimal controllers to explain DPF. We define and analyze DPF in the optimal control context, revisiting several control problems (state feedback, full control, and output feedback) to explore conditions that necessitate DPF. We find that even small deviations from the unconstrained state feedback problem (e.g. incomplete sensing, communication delay) necessitate DPF in the optimal controller. We also discuss parallels between controller structure and observations from neuroscience. In particular, the system level (SLS) controller displays DPF patterns compatible with predictive coding theory and easily accommodates signaling restrictions (e.g. delay) typical to neurons, making it a candidate for use in sensorimotor modeling.

\end{abstract}
	\section{INTRODUCTION} \label{sec:introduction}

\begin{figure}[h]
	\centering
	\includegraphics[width=6cm]{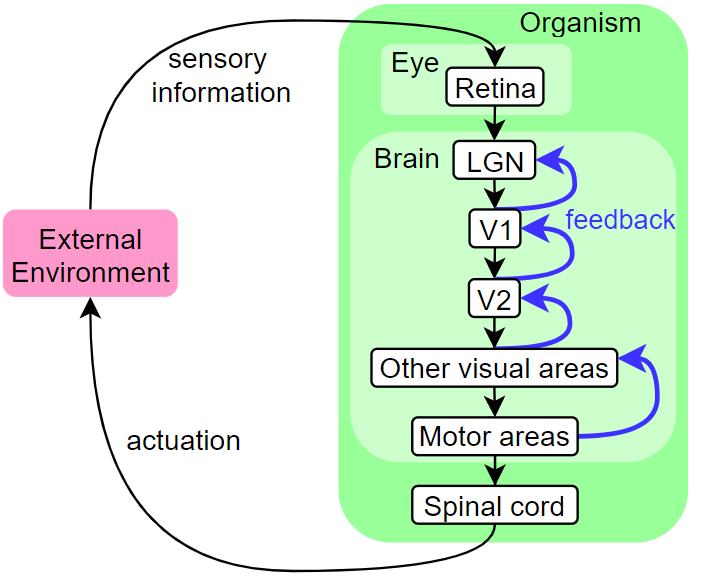}
	\caption{Simplified schematic of the human sensorimotor pathway. Blue arrows represent feedback, i.e. connectivity that defies the primary direction of information flow. There is additionally lots of lateral/loop connections within each box (e.g. V1, V2) as well, not shown in the diagram.} \label{fig:organism_fdbk}
\end{figure}

The primary visual pathway entails connections that propagate information from the retina in the eye to the lateral geniculate nucleus (LGN), then to the primary visual area (V1) in the cortex, secondary visual area (V2), and so on. However, massive amounts of connections in the reverse direction (i.e. descending predictive feedback or DPF, internal feedback, reciprocal connections) are also observed, as shown in Fig. \ref{fig:organism_fdbk}. This is a well-documented but poorly understood architectural feature \cite{Felleman1991, Callaway2004, Muckli2013}, and understanding the purpose of the mechanism is invaluable to understanding overall circuit function in the visual system. 

Current explanations surrounding visual DPF (e.g. modulation and memory processes, gain control, predictive coding, recurrent neural networks \cite{Muckli2013, Nayebi2018}) largely treat the visual system as an isolated module. However, we propose to view feedback from a sensorimotor perspective rather than a purely visual one. We are motivated by recent findings showing that non-visual-related neuronal activity in the visual areas -- which are substantial and thus far unexplained -- are dominated by body movements \cite{Stringer2019, Musall2019}. This suggests that motor activities play a significant role in neuronal activity in visual areas, and it is beneficial to analyze feedback in a sensorimotor context instead of a purely visual context. Optimal control models are an ideal candidate for this analysis as they have been widely successful in explaining behavior-level observations in the sensorimotor domain \cite{Todorov2004, Franklin2011}.

In optimal control applications, connectivity patterns within the controller are generally unexamined. Controllers are most often implemented on a computer or microcontroller, where connectivity merely translates to arithmetic operators in code. The exceptions are distributed control and delay-tolerant control, where connectivity patterns must work with communication constraints. To begin comparing optimal controller connectivity with neuronal connectivity, we must first clearly establish terminology surrounding \textit{feedback}; this is done in Section \ref{sec:terminology}. We provide an analysis of DPF within several canonical controllers in Section \ref{sec:int_fdbk_ctrl}. Notably, any deviation from the unconstrained state feedback problem (e.g. incomplete sensing, communication delay) necessitates DPF in the optimal controller. In \ref{sec:int_fdbk_sls}, we analyze DPF in a system level (SLS) controller \cite{Anderson2019} and note its compatibility with the popular predictive coding framework and other neuronal features. SLS easily accommodates signaling delay constraints, which are ubiquitous in neuronal systems, and is additionally scalable to large systems and thus suitable for producing higher-resolution models of the sensorimotor system. We discuss points of compatibility between controller structure and observations from neuroscience for all controllers analyzed.

Internal feedback is not a phenomenon unique to the visual system; it is prevalent across neuronal systems, and has been observed in the somatosensory and motor cortices \cite{Felleman1991}, auditory system \cite{Suga2008}, and other structures in the brain. Though our work focuses on the visual and motor systems, the analysis is sufficiently general to be of potential use in other domains as well. 

	\section{TERMINOLOGY} \label{sec:terminology}
We clarify the meaning of \textit{feedback} in the sensorimotor and control systems. In general, this term is highly overloaded, and the direction and scale of feedback may be unclear to researchers from different fields. We form our definitions of feedback using basic controls concepts, and show that it can be made consistent with notions of feedback in sensorimotor literature.

\begin{figure}
	\centering
	\includegraphics[width=6cm]{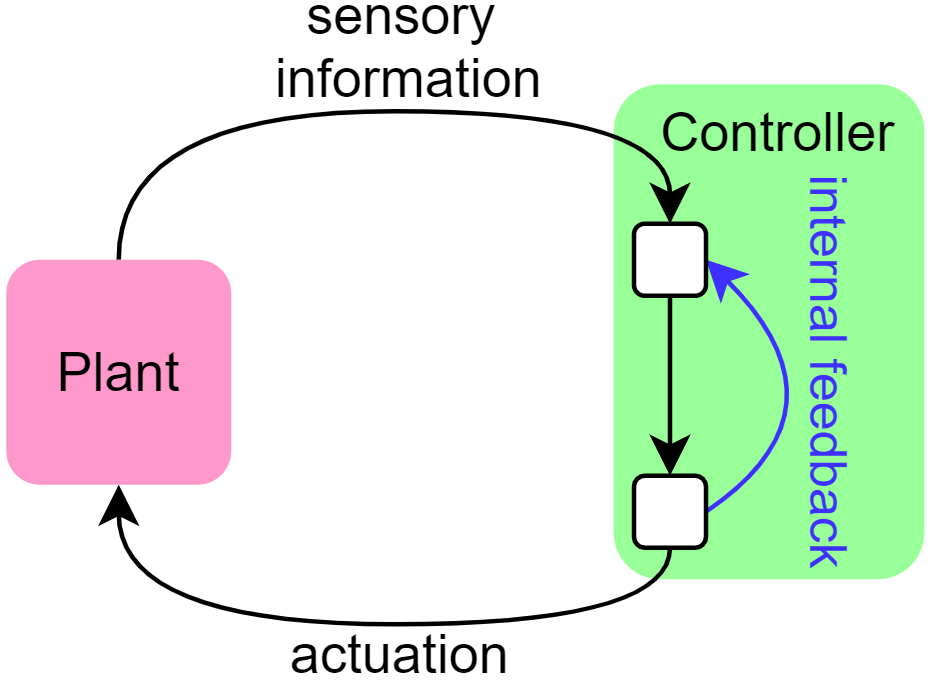}
	\caption{Block diagram of plant and controller. Internal feedback refers to information flow \textit{within} the controller, from the actuation back toward the sensory input. External feedback is represented by sensory information flow from the plant to the controller. This is the simplest diagram that contains internal feedback.} \label{fig:feedback_dir}
\end{figure}

\begin{definition}
	\textit{External feedback} refers to information flow from the plant to the controller.
\end{definition}

In the controls context, \textit{plant} and \textit{controller} are fairly obvious terms. In open-loop control, the controller acts on the plant in the absence of external feedback. In the sensorimotor setting, the external environment may be regarded as the plant, while the organism may be regarded as the controller. Here, external feedback can be thought of as information an organism receives from the external environment via its senses. This is shown by the `sensory information' arrow in Fig. \ref{fig:organism_fdbk}. We often refer to the combination of sensory information from plant to controller and actuation from controller to plant as a \textit{feedback loop} or a \textit{closed loop}.

In controls, feedback and external feedback are largely synonymous. External feedback is shown by the `sensory information' arrow in Fig. \ref{fig:feedback_dir}.  In closed-loop control, the controller makes use of this information to generate control signals. For example, in a system with state $x$, a control $u$ may be generated via state feedback by $u = Kx$ where $K$ is some matrix representing the controller.

\begin{definition}
	\textit{Internal feedback} refers to information flow within the controller (or estimator), traveling toward the sensor input.
\end{definition}

We note that a standalone estimator can still have internal feedback which travels from the estimator output toward the sensory input.

In the sensorimotor context, internal feedback refers to neuronal connectivity within the organism. In the visual setting, internal feedback refers to neural connections toward the sensor (i.e. eye); this is shown by the blue arrows in Fig. \ref{fig:organism_fdbk}. These connections are generally referred to as \textit{feedback} in the literature because they run opposite to the \textit{forward} (i.e. primary) direction of the visual pathway (i.e. retina to LGN to V1 to V2, and so on). They are also sometimes referred to as \textit{reciprocal} or \textit{bidirectional} connections.

In addition to forward and feedback connections, the brain also has \textit{lateral} connections; these are connections within components (e.g. V1 to V1), or connections across areas that are not definitively thought of as sequential relative to one another. Though we focus on forward and feedback connections, some ideas in this paper apply to lateral connections as well.

Internal feedback, as seen in the visual system, may alternatively be described as \textit{descending predictive feedback (DPF)}. In neuroscience, \textit{descending} generally refers to the direction away from cortex toward the periphery. In the visual system, descending feedback connects higher processing areas (e.g. V2) to lower processing areas (e.g. V1). Note that internal feedback is not always descending; in the motor system, feedback propagates from the spinal cord to the brain in an ascending direction. Feedback in the visual system is additionally described as \textit{predictive} as it is often assumed to serve a prediction-related role, such as predicting incoming sensory input based on internal models. 

Internal feedback in a simplified controller is shown by the blue arrow in Fig. \ref{fig:feedback_dir}. This is generally unexplored in technological applications, as internal feedback and controller connectivity are usually unimportant when said controller is implemented on a computer. However, understanding internal feedback in controllers allows us to draw interesting parallels with internal feedback in the visual system, as we will show in the following sections.

For the remainder of this paper we will use `internal feedback' and `DPF' interchangeably. Strictly speaking, DPF is a subset of internal feedback; however, our analysis of controller structure suggests that most internal feedback indeed serves a prediction-related purpose, coinciding with prevailing theories about vision. We tentatively adopt this terminology and hope that future work will provide sharper results that further refine definitions surrounding internal feedback in both controllers and sensorimotor systems.
 
	\section{INTERNAL FEEDBACK IN CANONICAL CONTROLLERS} \label{sec:int_fdbk_ctrl}

Using simple canonical controllers and examining connectivity patterns within the controller, we observe that though the basic state feedback problem requires no internal feedback, any amount of incomplete sensing or communication delay in the problem will produce a controller containing internal feedback. Generally speaking, internal feedback seems to serve the purpose of passing information about the control action $u$ toward some internal estimated quantity (e.g. state estimate) to improve estimation of said quantity. Often, the estimated quantity is indirectly affected by $u$; however, it is preferable to have direct knowledge of $u$ rather than infer $u$ through the output $y$. Improved estimation improves the quality of future control actions and performance.

The qualitative ideas above are not rigorous, and our overall analysis of internal feedback is not comprehensive. Internal feedback within controllers is generally unexamined; we aim to change this by presenting canonical examples to lay the foundation for further theoretical work, and by motivating our analysis with connections to neuroscience. Since optimal control models are powerful at explaining sensorimotor behaviors \cite{Todorov2004} and display prevalent DPF, we believe that they can explain some portion of DPF in sensorimotor structures, including the visual system.

\subsection{Scalar Controllers}
We use the following plant, where $x$ is the external state, $u$ is the actuation exerted upon the external state by the controller, and $w$ is disturbance on the external state. $w$ is assumed to be white noise with variance $\sigma_w^2$. 

\begin{equation}
	x(t+1) = ax(t) + u(t) + w(t)
\end{equation}

We sense the external state plus some white noise $v(t)$ with variance $\sigma_v^2$.

\begin{equation}
	y(t) = x(t) + v(t)
\end{equation}

When $\sigma_v^2=0$, this can be viewed as a state feedback (SF) or full control (FC)\footnote{Full control is applicable to any problem with full actuation, and is the dual problem of state feedback (which is applicable to any problem with perfect sensing). The full control gain is equal to the observer gain for output feedback, just as the state feedback gain is equal to the controller gain for output feedback.} problem. In both cases, the optimal gain $k$ (or $l$ for FC) is obtained by solving a discrete algebraic Riccati equation (DARE), and the optimal controller is $u(t) = -ky(t)$ ($u(t) = -ly(t)$ for FC). No internal feedback is present.

For nonzero sensing noise, SF no longer applies, so we must use the FC controller. As before, no internal feedback is present. The two simple gain controllers are shown in Fig. \ref{fig:sf_fc_of}.

\subsection{Estimator Dynamics and Output Feedback}
\begin{figure}
	\centering
	\includegraphics[width=8.5cm]{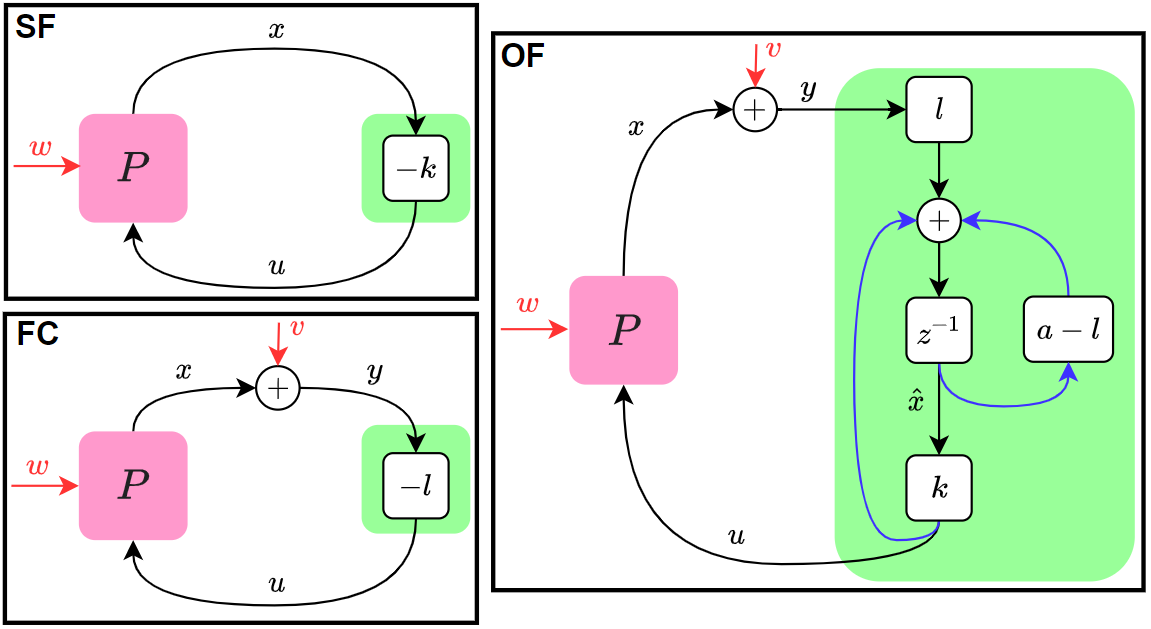}
	\caption{Block diagrams of state feedback (SF), full control (FC), and output feedback (OF) controllers. The green rectangle is the controller. DPF in the OF controller is denoted by blue arrows.} \label{fig:sf_fc_of}
\end{figure}

For nonzero noise, we can also view this problem as an output feedback (OF) problem. The optimal solution combines an optimal controller $k$ with an optimal estimator $l$. The presence of the optimal estimator necessitates internal dynamics, which contain DPF. The controller is implemented as

\begin{subequations}
	\begin{equation}
		\hat{x}(t+1) = (a-l)\hat{x}(t) + u(t) + ly(t)
	\end{equation}
	\begin{equation}
		u(t) = -k\hat{x}(t)
	\end{equation}
\end{subequations}

\noindent where $k$ and $l$ are the gains for the SF and FC problem, respectively. The structure of the optimal SF, FC, and OF controllers are shown in Fig. \ref{fig:sf_fc_of}.

Two internal feedback connections are observed: one carrying the $u$ signal, and one passing through the $a-l$ block.

\begin{subequations} \label{eqn:output_fdbk}
	\begin{equation}
	\hat{x}(t+1) = (A-LC)\hat{x}(t) + Bu(t) + Ly(t)
	\end{equation}
	\begin{equation}
	u(t) = -K\hat{x}(t)
	\end{equation}
\end{subequations}

More generally, in the generic optimal OF controller for plant $(A, B, C)$ shown in \ref{eqn:output_fdbk}, the gains $B$ and $(A-LC)$ constitute DPF. For nontrivial actuation matrix $B$, an OF controller will \textit{always} have DPF due to its estimator dynamics. This is true in the continuous time case as well, in which $B$ and $(A-LC)$ again constitute DPF. 

In both discrete and continuous time, OF contains two DPF terms whose functions we list below. Note that even without a closed-loop control system, a solitary estimator (i.e. Kalman filter) would still contain internal feedback in the form of estimator dynamics.
\begin{enumerate}
	\item Updating the estimator with motor plans via $B$
	\item Estimator dynamics: accounting for anticipated plant dynamics and estimator error via $L$ and an internal copy of $A$ and $C$ (i.e. $A-LC$)
\end{enumerate}

Note: we include OF for the purpose of demonstrating estimator dynamics on a simple system. For our system, OF is not strictly required since FC suffices. In general, OF will be necessary when we have neither perfect sensing nor full actuation, as is true in the human sensorimotor system. Additionally, though we focus on $\mathcal{H}_2$ control in this work, internal feedback is present in $\mathcal{H}_{\infty}$ control as well. In the OF case, the $\mathcal{H}_{\infty}$ controller includes the two above-mentioned terms plus an additional feedback term of predicted worst-case disturbance.

\subsubsection{Efference copy}
In the sensorimotor context, updating an estimator with motor plans is a familiar idea  \cite{Todorov2004}. In neuroscience terminology, a copy of motor output -- termed the \textit{efference copy} or \textit{corollary discharge} -- is hypothesized to be sent to an internal estimator to compensate for sensory delay. We note that even in a delay-free OF controller, internal feedback similar to the efference copy is present. In the brain, the efference copy along with other sensory inputs contributes to some type of Bayesian estimation \cite{Franklin2011}. The location of state estimation within the brain is an open question; the posterior parietal cortex and cerebellum \cite{Desmurget2000, Shadmehr2008} are two candidates.

\subsubsection{Kalman filter and forward models}
A standalone Kalman filter contains internal feedback in the form of estimator dynamics. Broadly speaking, the visual cortex parses pixels detected by the retina into meaningful entities (e.g. objects, organisms, and environments), which requires some form of estimation and filtering akin to a Kalman filter. Some amount of visual internal feedback is likely explained by this, though findings of motor activity in visual cortex \cite{Stringer2019, Musall2019} suggest that some amount of visual internal feedback may also be related to motor plans. Additionally, this explanation is incomplete in the sensorimotor domain as it does not incorporate any notion of signaling delay, which is dominant in neuronal systems. The concept of accounting for plant dynamics is also related to that of a forward model (i.e. an internal copy of $A$ and $B$), which uses the current state and motor output to guess future estates. In a macroscopic sense, human sensorimotor behavior appears to incorporate such a forward model, yet it is difficult to conclusively demonstrate the existence of the model within the brain \cite{Franklin2011}. Lastly, we note that the $(A-LC)$ estimator dynamics in the OF controller can be interpreted as lateral connections.

\subsubsection{Size of forward and feedback signals}
In the OF controller, the feedback signal that carries motor plans has dimension equal to input size, while the signal that accounts for plant dynamics has dimension equal to state size. Together, the information carried by internal feedback in OF has higher dimension than the information carried by the forward path, which has dimension equal to state size. This qualitatively agrees with findings that the feedback path contains more axons than the feedforward path, e.g. V1 receives 10 times more axons from V2 (i.e. internal feedback) than from LGN (i.e. feedforward) \cite{Budd1998}.

\subsection{Delayed Communication and Full Control}
\begin{figure}
	\centering
	\includegraphics[width=7cm]{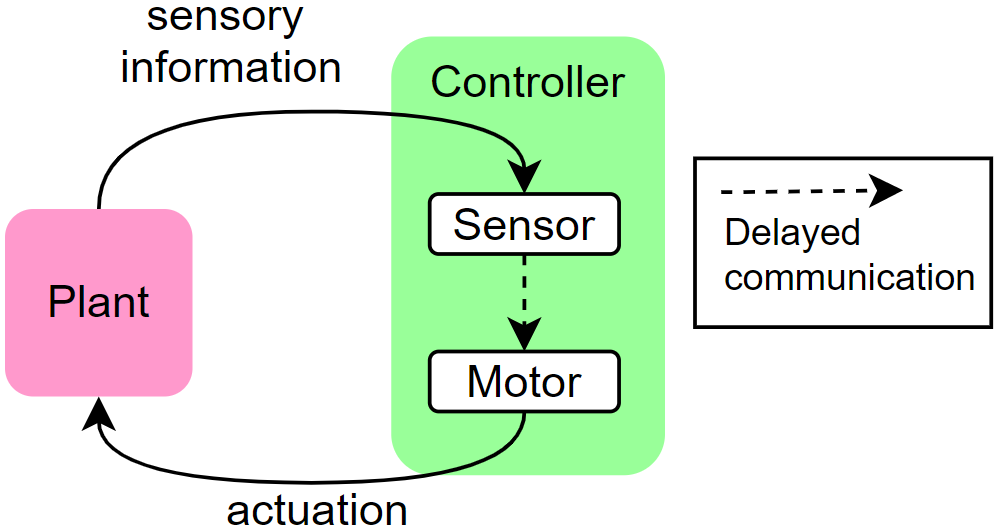}
	\caption{Diagram of controller with delayed communication between sensor and motor regions. No internal feedback is depicted; internal feedback, if present, would run in the motor-to-sensor direction.} \label{fig:sensor_motor}
\end{figure}

We now consider the case where communication within the controller is delayed. In particular, we assume delay between the region of the controller that senses (`sensor') and the region of the controller that computes the control action (`motor'); this is shown in Fig. \ref{fig:sensor_motor}. This type of delay plays a big role in the sensorimotor system and brain in general, where delayed signaling abounds due to the physical limitations of spiking neurons.

We rename the external state $x_1$. Now, consider a delay state $x_2$, which lags the true state by one time step:

\begin{subequations}
	\begin{equation}
		x_2(t+1) = x_1(t)
	\end{equation}
	\begin{equation}
		y(t) = x_2(t)
	\end{equation}
\end{subequations}

The sensor region senses this true state in real time and communicates it to the motor region using a delayed wire. Thus, the motor region accesses the delay state instead of the true state. We show that in the presence of this delay, the optimal controller must have DPF. First, we write the system in standard state space form:

\begin{subequations}
\begin{equation} \label{eq:standard_statespace}
	x(t+1) = Ax(t) + Bu(t) + w(t)
\end{equation}
\begin{equation}
	y(t) = Cx(t) + v(t)
\end{equation}
\end{subequations}

\noindent with state $x = \begin{bmatrix} x_1 & x_2 \end{bmatrix}^\top$

\noindent and matrices
\begin{equation} \label{eq:state_space}
	A = \begin{bmatrix} a & 0 \\ 1 & 0 \end{bmatrix} \quad B = \begin{bmatrix} 1 \\ 0 \end{bmatrix} \quad C = \begin{bmatrix} 0 & 1 \end{bmatrix}
\end{equation}

This problem can be solved using OF, which will yield DPF as per previous discussion. We may also reframe it as an FC problem.

Initially, this problem appears incompatible with FC since we do not have full actuation, i.e. $B \neq I$. However, since our problem formulation treats $x_2$ as an \textit{internal} state, i.e. as part of the controller, we may add `actuation' to $x_2$ to make $B = I$. This translates to adding an extra wire in the controller, which is acceptable since we are free to design the controller and an extra communication wire is cheap. Note that we are only able to do this because we defined $x_2$ to be an internal state; if it were an external state, it would be costly or physically impossible to force additional actuation on it, and we would be unable to make $B = I$.

In the FC formulation, the delay state becomes
\begin{equation}
	x_2(t+1) = x_1(t) + u_2(t)
\end{equation}

\noindent where $u_2$ is an \textit{internal} control signal. We define the control $u$ as \noindent with state $u = \begin{bmatrix} u_1 & u_2 \end{bmatrix}^\top$

Although both $u_1$ and $u_2$ are control signals in the standard sense, they have vastly different interpretation and cost; $u_1$ represents costly physical actuation, while $u_2$ represents a cheap communication wire. Additionally, the noise $v$ on the output $y$ can now be interpreted as noise internal to the controller, i.e. due to noisy signaling. The optimal controller can be written as
\begin{equation}
\begin{bmatrix} u_1 \\ u_2 \end{bmatrix} = -\begin{bmatrix} l_1 \\ l_2 \end{bmatrix} y
\end{equation}

\begin{figure}
	\centering
	\includegraphics[width=7cm]{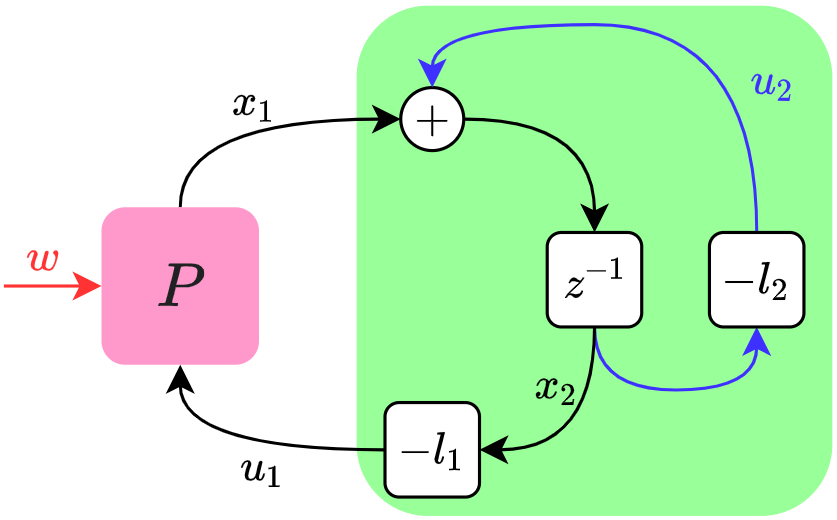}
	\caption{FC controller for system with a delay of one time-step, represented by the $z^{-1}$ block. All signals are scalar. $k_i$ represents the $i$th element of the controller gain, i.e. $K = \begin{bmatrix}k_1 & k_2\end{bmatrix}$. For this system, analytical solutions show that $k_2$ is always zero. $l_i$ represents the $i$th element of the observer gain, i.e. $L = \begin{bmatrix}l_1 & l_2\end{bmatrix}^{\top}$. } \label{fig:delay_system}
\end{figure}

The FC controller is shown in Fig. \ref{fig:delay_system}. Since we consider $x_2$ and $u_2$ to be internal signals, the FC controller displays DPF, which is represented by the $l_2$ gain. We now examine two cases when this DPF is necessary, i.e. when $l_2 \neq 0$. 

First, DPF is necessary for optimal closed-loop performance. From the DARE, we obtain that

\begin{subequations}
\begin{equation}
	l_2 = \frac{p_2}{p_2+\sigma_v^2}, \quad l_1 = al_2
\end{equation}

\begin{equation}
	p_2 = \frac{a^2}{2}\left(-\beta \pm \sqrt{\beta^2 + \left( \frac{2}{a}\sigma_w\sigma_v\right)^2} \right)
\end{equation}

\begin{equation}
	\beta = (\frac{1}{a} - a)\sigma_v^2 - \frac{1}{a}\sigma_w^2
\end{equation}
	
\end{subequations}

Here, $p_2$ represents the off-diagonal entries of the symmetric matrix $P$ which solves the DARE, chosen so that $P \succ 0$. We notice that $l_2=0$ implies $l_1=0$; when no DPF is required, no control action is required, i.e. $\mu=0$ is the optimal choice. In this case, the loop is effectively open. Optimal performance in this closed-loop system appears to require DPF. In this controller, $x_2$ roughly represents an estimate of current state, adjusted by estimated dynamic propagation and noise factors via $u_2$. $x_2$ is present due to delay; without delay, the FC controller would have no DPF as discussed in earlier sections.

\begin{figure}
	\centering
	\includegraphics[width=7cm]{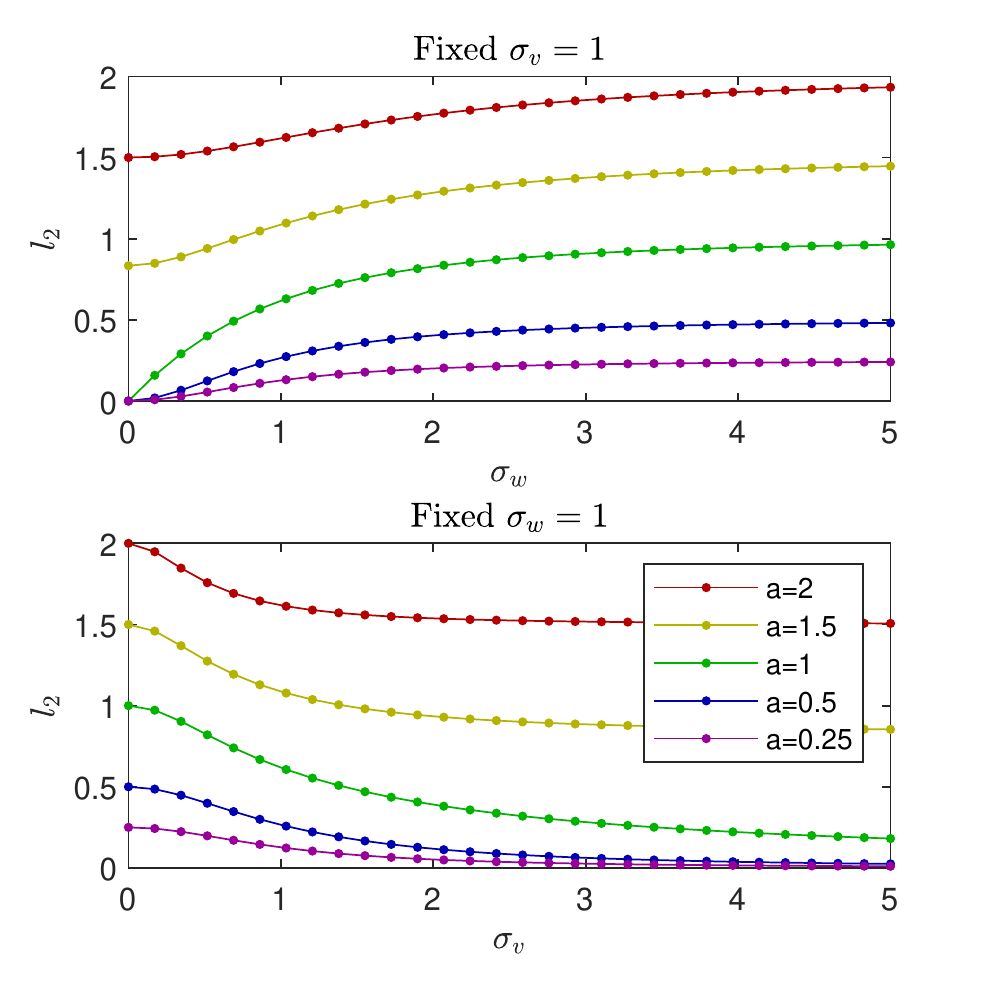}
	\caption{DPF gain $l_2$ for various stability values vs. disturbance variance ($sig_w$) and internal noise ($sig_v$). Gain magnitude increases with disturbance variance, and decreases with stability and noise.} \label{fig:l2_plot}
\end{figure}

We additionally inspect the behavior of $l_2$ for different values of $a$, $\sigma_w$, and $\sigma_v$ in Fig. \ref{fig:l2_plot}. Here, we restrict disturbance $w$ to only act on external state $x_1$. We see that $l_2$ is only zero when $a=0$ or $\sigma_w=0$ and $|a| < 1$; the optimal controller always employs DPF for unstable systems, and for stable systems when external disturbance is expected and $a \neq 0$.

DPF is not only crucial for optimal performance but also for stabilizing an unstable open-loop system. For this simple case with a delay of one time-step, a nonzero $l_2$ is necessary to stabilize unstable systems with $|a| \geq 2$. This analysis extends to systems with larger delays. For example, instead of having a single delay state, we can use a net delay of $T_d$ time-steps and add additional delay states

\begin{subequations}
	\begin{equation}
		x_i(t+1) = x_{i-1}(t) + u_i(t), i = 2 \ldots T_d+1
	\end{equation}
	\begin{equation}
		y(t) = x_{T_d+1}(t)
	\end{equation}
\end{subequations}

\begin{figure}
	\centering
	\includegraphics[width=5cm]{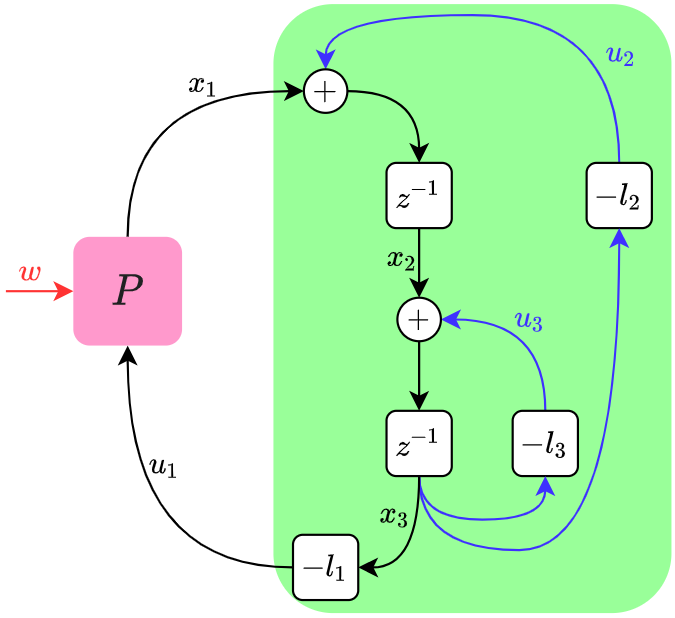}
	\caption{FC controller for system with a net sensing delay of two timesteps, i.e. $T_d=2$. $l_i$ represents the $i$th element of the observer gain, i.e. $L = \begin{bmatrix}l_1 & l_2 & l_3\end{bmatrix}^{\top}$.} \label{fig:two_delay}
\end{figure}

\begin{figure}
	\centering
	\includegraphics[width=7cm]{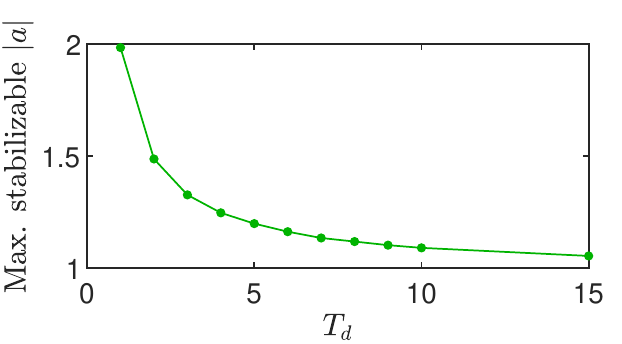}
	\caption{Maximum stabilizable $|a|$ for various net delays $T_d$, without DPF. With DPF, the system can be stabilized for all $a$, i.e. the maximum stabilizable $|a|$ is unbounded.} \label{fig:stab_without_intfdbk}
\end{figure}

A system with $T_d=2$ is shown in Fig. \ref{fig:two_delay}. For delayed systems, gains $l_i, i \geq 2$ represent DPF. It is generally true that for these higher delay systems, the optimal controller employs DPF, i.e. $\exists i, i \geq 2, l_i \neq 0$. Furthermore, systems with larger net delays also require DPF  to stabilize some unstable open-loop systems. We plot the maximum $|a|$ that a controller with no DPF can stabilize (i.e. $l_i=0, i=2 \ldots T_d$) in Fig. \ref{fig:stab_without_intfdbk}. As delay becomes larger, the set of systems that are stabilizable without DPF becomes smaller. Even at $T_d=2$, we need DPF to stabilize systems with $|a| \geq 1.5$. Note that with DPF, we can stabilize any plant because for any delay, $(C, A)$ will be observable and thus we can arbitrarily place the closed-loop eigenvalues.

\subsubsection{Neuronal signaling delay}
\begin{figure}
	\centering
	\includegraphics[width=7cm]{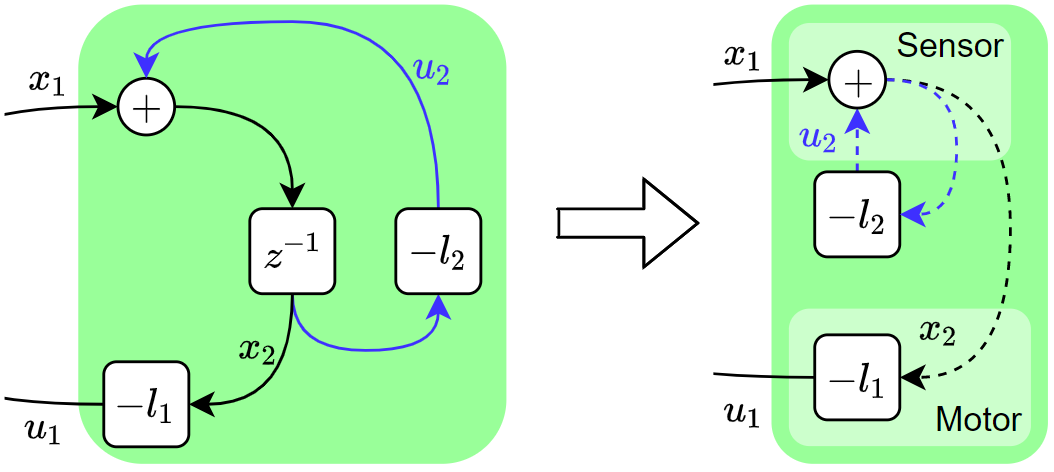}
	\caption{Two equivalent controllers. The controller on the right has been partitioned into two regions (`sensor' and `motor') with inter-regional delay. Dashed lines indicate delayed communication. The black dashed line corresponds to 1 timestep ($z^{-1}$) of delay. The blue dashed lines correspond to 1 timestep of delay overall.} \label{fig:fc_realization}
\end{figure}

Delayed signaling is abundant in the nervous system. Neurons can send faster signals by increasing axon diameter but this incurs significant metabolic cost, especially over long distances \cite{Sterling2015}; for this reason, signaling between regions that are farther apart tends to incur more delay. Even over shorter distances, signals may be sent slowly to conserve energy. In the simple sensor/motor system in Fig. \ref{fig:sensor_motor}, we assume that intra-region communication (i.e. connections within `sensor' or `motor') has negligible delay, while inter-region communication (i.e. communications between sensor and motor units) is subject to delay. We redraw the one-step delay FC-controller from Fig. \ref{fig:delay_system} to reflect these delays in Fig. \ref{fig:fc_realization}. In this diagram, the location of the DPF term $l_2$ is ambiguous; it may reside somewhere between the motor and sensor regions or be a part of slow communication within the sensor region. The latter explanation has some parallels with DPF -- including lateral connections -- in the visual system, including findings that feedback paths use slower receptors than forward paths \cite{Attwell2005, Self2012}.

\subsubsection{Size of forward and feedback signals}
In the FC controller, the DPF signal has dimension equal to state dimension multiplied by steps of communication delay. The forward signal has the dimension of the observed quantity, which is generally upper-bounded by state dimension. As with the OF controller, the feedback signal has higher dimension than the forward signal, mirroring findings in neuroscience. For high delay, the FC controller predicts a bigger difference between feedforward and feedback signal dimension compared to the OF controller. Nonetheless, for the same problem, the FC controller is algebraically simpler and requires fewer signaling wires than the OF controller.

Unlike the OF controller, the FC controller does not explicitly estimate the state. Instead, the FC controller implicitly performs estimation and accounts for delays, plant dynamics, and motor plans. Thus, the discussion on efference copy and forward models from the previous section can be applied to the FC controller as well.

	\section{INTERNAL FEEDBACK IN SYSTEM-LEVEL CONTROLLERS} \label{sec:int_fdbk_sls}

So far, we have analyzed canonical controllers in simple settings to extrapolate a principle surrounding the presence of internal feedback. The FC controller for a delayed system has basic structural similarities to visual DPF, but requires some hand-crafting. To formulate a problem with internal delay as an FC problem, we must incorporate internal delay states into the plant. For a simple scalar system, this is easy, but this formulation becomes cumbersome if we want to rapidly experiment with varying delays on a more complex system. 

Delayed signaling plays an important role in the visual system and neuronal systems in general, and sets them apart from standard controllers; most control applications use fast electronics and do not face significant signaling delays. Neuronal communication delay is used in \cite{Nakahira2015, Nakahira2018, Liu2019} to motivate a layered controller architecture in the sensorimotor system, though detailed controller connectivity is not examined in these works. We eventually aim to incorporate analysis of DPF into the layered framework presented in the abovementioned works. To begin doing so, we must first incorporate delayed communication into a controller in a scalable way; we do so using the recently proposed system-level controller. 

The system level controller is part of the \textit{System Level Synthesis} (SLS) framework \cite{Anderson2019}. We will refer to it as the SLS controller. This controller is promising for producing sensorimotor models because it easily accommodates signaling restrictions (e.g. delay) typical to neurons and scales well. Like the OF and FC controllers, the SLS controller has several points of qualitative agreement with existing theories and findings in neuroscience.

The standard SLS controller for both SF and OF utilize DPF. Given previous discussion, DPF in OF should not be surprising. We focus on the SF-SLS controller, emphasizing controller implementation rather than controller synthesis. For a full setup of the synthesis problem and scalability analysis, we refer the reader to $\S4$ of \cite{Anderson2019}, where the derivation and benefits of the SLS approach are thoroughly discussed. 

First, we write equation (\ref{eq:standard_statespace}) in the frequency domain as

\begin{equation} \label{eq:freq_statespace}
z\x = A\x + B\u + \w
\end{equation}

\noindent with controller
\begin{equation} \label{eq:freq_ctrller}
\u = \K\x
\end{equation}

\noindent where $\K$ is some transfer matrix. Instead of directly searching over possible controllers $\K$, the SLS approach searches over $\left\{\R, \M\right\}$, the closed-loop responses of state and control input to disturbances, i.e.

\begin{equation}
\begin{bmatrix} \x \\ \u \end{bmatrix} = \begin{bmatrix} \R \\ \M \end{bmatrix} \mathbf{w}
\end{equation}	

To ensure that a controller $\K$ exists to achieve these closed-loop responses, we enforce
\begin{equation} \label{eq:sls_constraint}
\begin{bmatrix} zI-A & -B \end{bmatrix} \begin{bmatrix} \R \\ \M \end{bmatrix} = I
\end{equation}

\begin{figure}
	\centering
	\includegraphics[width=6cm]{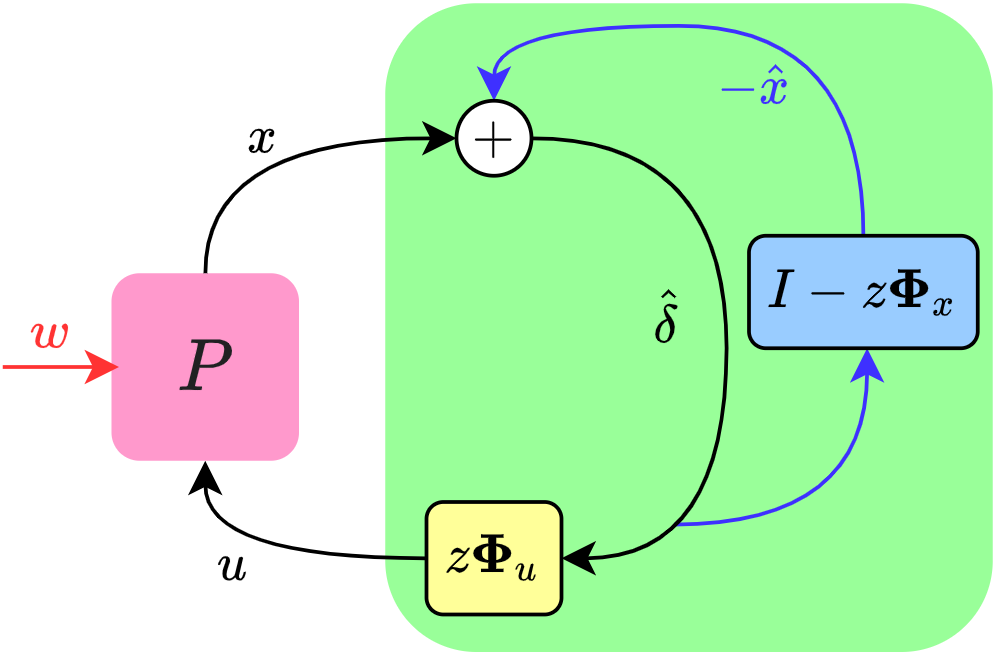}
	\caption{Block diagram for state feedback SLS controller. Signals $x$ and $u$ are vectors and the blocks represent transfer matrices instead of scalar gains. The blue $I-z\R$ block represents DPF.} \label{fig:sls_basic}
\end{figure}

The controller is then implemented using $\left\{\R, \M\right\}$ as shown in Fig. \ref{fig:sls_basic}. DPF in the controller is central to the SLS parametrization, as it enables formulation of previously difficult constraints (e.g. delayed and local communication) as affine constraints on  $\left\{\R, \M\right\}$. This structure also gives easy-to-interpret signals; $\hat{x}$ represents a prediction of the next expected state, which is subtracted from the incoming sensed state; $\hat{\delta}$ is the resulting difference, representing estimated disturbance which was not predicted by DPF. The forward transfer matrix $z\M$ then acts upon this difference signal instead of the full state signal. We note that in a case with no delay, the SLS structure is not required and a scalar-gain SF controller with no DPF suffices. However, as soon as any communication constraints come into play, the SLS structure with DPF becomes necessary.

As mentioned in the previous section, a key feature of neuronal signaling is communication delay. The SLS formulation easily accommodates these through sparsity constraints on elements of  $\left\{\R, \M\right\}$. To demonstrate this, we first write out the transfer matrices in terms of their spectral components:

\begin{subequations}
	\begin{equation}
	\R(z) = z^{-1}I + z^{-2}\Phi_x^{(2)} + z^{-3}\Phi_x^{(3)}
	\end{equation}
	\begin{equation}
	\M(z) = z^{-1}\Phi_u^{(1)} + z^{-2}\Phi_u^{(2)} + z^{-3}\Phi_u^{(3)}
	\end{equation}
\end{subequations}

\begin{figure}
	\centering
	\includegraphics[width=5cm]{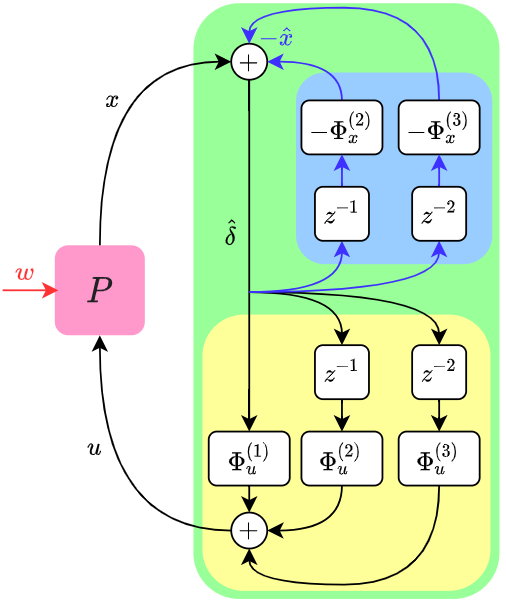}
	\caption{Expanded block diagram for state feedback SLS controller. All blocks are scalar gains. The blue and yellow boxes correspond to the $I-z\R$ and $z\M$ blocks from Fig. \ref{fig:sls_basic}, respectively. The $z^{-1}$ block indicates one step of delay; $z^{-2}$ indicates two steps of delay.} \label{fig:sls_expanded}
\end{figure}

For simplicity, we restrict $\left\{\R, \M\right\}$ to three spectral components each. The constraint (\ref{eq:sls_constraint}) enforces that $\Phi_x^{(1)} = I$. In general, $\R$ and $\M$ are strictly proper finite impulse response transfer matrices with $T$ spectral components for some $T \geq 1$, i.e.
\begin{subequations}
	\begin{equation}
	\R(z) = \sum_{k=1}^{T} \Phi_x^{(k)}z^{-k}
	\end{equation}
	\begin{equation}
	\M(z) = \sum_{k=1}^{T} \Phi_u^{(k)}z^{-k}
	\end{equation}
\end{subequations}

\noindent where $\Phi_x^{(1)} = I$. The controller is shown in Fig. \ref{fig:sls_expanded}. 

$\Phi_u^{(1)}$ represents non-delayed forward information propagation (i.e. from sensor toward motor). $\Phi_u^{(2)}$ and  $\Phi_u^{(3)}$ represent forward information propagated with one and two timesteps of delay, respectively. Similarly,  $\Phi_x^{(2)}$ and $\Phi_x^{(3)}$ represent feedback information propagated with one and two timesteps of delay, respectively. Note that the SLS formulation inherently does not include a non-delayed feedback propagation term; this term is canceled out by $I$ in the $I-z\R$ block.

To enforce some minimum delay, we only have to assert that the appropriate gains are zero. Additionally, we can accommodate arbitrary amounts of delay with additional spectral components. This is in contrast with the FC model of delay, in which delays must be hand-crafted into the plant model.

\subsubsection{Predictive coding}
The SF SLS controller structure shown in Fig. \ref{fig:sls_basic} resembles predictive coding theory in neuroscience \cite{Huang2011}, which has strong experimental support in the visual areas. Predictive coding hypothesizes that higher processing stages (e.g. V2) learn to predict input to lower processing stages (e.g. V1) and communicate this prediction via DPF; unpredicted quantities (e.g. errors) are communicated via forward connections. This corroborates findings that uninstructed movements account for much neural activity in visual cortex \cite{Stringer2019, Musall2019}. Additionally, this looks like what the SLS controller does. The $\R$ term can be interpreted as an internal predictive model in this scenario -- recent work combining SLS and learning indeed treats $\R$ (and $\M$) as a learned quantity \cite{Xue2020}, which again coincides with predictive coding theory. While predictive coding is motivated by information efficiency, the SLS controller is motivated by optimal performance under communication constraints. The convergence of the two suggest that the optimal control explanation is compatible with, rather than at odds with, prevailing notions of DPF in the neuroscience community. We remark that the resemblance to predictive coding theory also holds in the OF-SLS controller, which is more structurally complex and not depicted here. 

\subsubsection{Size of forward and feedback signals}
In SF-SLS, the internal forward signal has dimension equal to state dimension. The feedback signal has dimension roughly equal to state dimension multiplied by steps of communication delay. This dimension may be reduced by limiting the number of paths with different delays. Like the OF and FC controllers, the SF-SLS controller has a feedback signal with higher dimension than the forward signal, reflecting neuroscience results. The feedback signal dimension is increased when the SF-SLS controller is implemented in a distributed manner, i.e. each subsystem in the system implements a local controller that only uses local information. In this case, each state's information is stored by all local controllers which access it, creating some duplication of information.

\subsubsection{Sensorimotor models}
The SLS controller is promising for producing sensorimotor models as it easily accommodates signaling delay. As previously discussed, signaling delay is prevalent in neuronal systems. There is ample evidence of diverse conduction speeds within the visual system \cite{Angelucci2003, Perge2009}; SLS can accommodate this as it allows for multiple pathways with different delays. Originally conceived as a distributed controls tool, SLS is also scalable to large systems and thus suitable to model the sensorimotor system in higher resolution.

	\section{CONCLUSION} \label{sec:conclusion}

We defined the notion of internal feedback and provided a controls-based analysis of its presence and necessity using a variety of simple controllers. Each controller has some agreement with concepts and observations in neuroscience, though the SLS controller appears to have the most promise for producing sensorimotor models due to its scalability and ability to easily accommodate arbitrary communication delays. This paper provides broad theoretical analysis on internal feedback that hopefully complements data-driven techniques in neuroscience and invokes discussion and cooperation between control theorists and neuroscientists. Future work will aim to produce more complex and realistic sensorimotor models using SLS to more quantitatively connect to observations in neuroscience.
	
	\bibliography{internal_feedback}

\begin{thebibliography}{10}
\providecommand{\url}[1]{#1}
\csname url@samestyle\endcsname
\providecommand{\newblock}{\relax}
\providecommand{\bibinfo}[2]{#2}
\providecommand{\BIBentrySTDinterwordspacing}{\spaceskip=0pt\relax}
\providecommand{\BIBentryALTinterwordstretchfactor}{4}
\providecommand{\BIBentryALTinterwordspacing}{\spaceskip=\fontdimen2\font plus
\BIBentryALTinterwordstretchfactor\fontdimen3\font minus
  \fontdimen4\font\relax}
\providecommand{\BIBforeignlanguage}[2]{{%
\expandafter\ifx\csname l@#1\endcsname\relax
\typeout{** WARNING: IEEEtran.bst: No hyphenation pattern has been}%
\typeout{** loaded for the language `#1'. Using the pattern for}%
\typeout{** the default language instead.}%
\else
\language=\csname l@#1\endcsname
\fi
#2}}
\providecommand{\BIBdecl}{\relax}
\BIBdecl

\bibitem{Felleman1991}
D.~J. Felleman and D.~C. {Van Essen}, ``{Distributed hierarchical processing in
  the primate cerebral cortex},'' \emph{Cerebral Cortex}, vol.~1, no.~1, pp.
  1--47, 1991.

\bibitem{Callaway2004}
E.~M. Callaway, ``{Feedforward, feedback and inhibitory connections in primate
  visual cortex},'' \emph{Neural Networks}, vol.~17, no. 5-6, pp. 625--632,
  2004.

\bibitem{Muckli2013}
L.~Muckli and L.~S. Petro, ``{Network interactions: non-geniculate input to
  V1},'' \emph{Current Opinion in Neurobiology}, vol.~23, no.~2, pp. 195--201,
  2013.

\bibitem{Nayebi2018}
A.~Nayebi, D.~Bear, J.~Kubilius, K.~Kar, S.~Ganguli, D.~Sussillo, J.~J.
  DiCarlo, and D.~L. Yamins, ``{Task-driven convolutional recurrent models of
  the visual system},'' in \emph{NeurIPS}, 2018, pp. 5295--5306.

\bibitem{Stringer2019}
C.~Stringer, M.~Pachitariu, N.~Steinmetz, C.~B. Reddy, M.~Carandini, and K.~D.
  Harris, ``{Spontaneous behaviors drive multidimensional, brainwide
  activity},'' \emph{Science}, vol. 364, 2019.

\bibitem{Musall2019}
S.~Musall, M.~T. Kaufman, A.~L. Juavinett, S.~Gluf, and A.~K. Churchland,
  ``{Single-trial neural dynamics are dominated by richly varied movements},''
  \emph{Nature Neuroscience}, vol.~22, no.~10, pp. 1677--1686, 2019.

\bibitem{Todorov2004}
E.~Todorov, ``{Optimality principles in sensorimotor control},'' \emph{Nature
  Neuroscience}, vol.~7, no.~9, pp. 907--915, 2004.

\bibitem{Franklin2011}
D.~W. Franklin and D.~M. Wolpert, ``{Computational mechanisms of sensorimotor
  control},'' \emph{Neuron}, vol.~72, no.~3, pp. 425--442, 2011.

\bibitem{Anderson2019}
J.~Anderson, J.~C. Doyle, S.~H. Low, and N.~Matni, ``{System level
  synthesis},'' \emph{Annual Reviews in Control}, vol.~47, pp. 364--393, 2019.

\bibitem{Suga2008}
N.~Suga, ``{Role of corticofugal feedback in hearing},'' \emph{Journal of
  Comparative Physiology A: Neuroethology, Sensory, Neural, and Behavioral
  Physiology}, vol. 194, no.~2, pp. 169--183, 2008.

\bibitem{Desmurget2000}
M.~Desmurget and S.~Grafton, ``{Forward modeling allows feedback control for
  fast reaching movements},'' \emph{Trends in Cognitive Sciences}, vol.~4,
  no.~11, pp. 423--431, 2000.

\bibitem{Shadmehr2008}
R.~Shadmehr and J.~W. Krakauer, ``{A computational neuroanatomy for motor
  control},'' \emph{Experimental Brain Research}, vol. 185, no.~3, pp.
  359--381, 2008.

\bibitem{Budd1998}
J.~M. Budd, ``{Extrastriate feedback to primary visual cortex in primates: A
  quantitative analysis of connectivity},'' \emph{Proceedings of the Royal
  Society B: Biological Sciences}, no. 1400, pp. 1037--1044, 1998.

\bibitem{Sterling2015}
P.~Sterling and S.~B. Laughlin, \emph{{Principles of neural design}}.\hskip 1em
  plus 0.5em minus 0.4em\relax MIT Press, 2015.

\bibitem{Attwell2005}
D.~Attwell and A.~Gibb, ``{Neuroenergetics and the kinetic design of excitatory
  synapses},'' \emph{Nature Reviews Neuroscience}, vol.~6, no.~11, pp.
  841--849, 2005.

\bibitem{Self2012}
M.~W. Self, R.~N. Kooijmans, H.~Sup{\`{e}}r, V.~A. Lamme, and P.~R. Roelfsema,
  ``{Different glutamate receptors convey feedforward and recurrent processing
  in macaque V1},'' \emph{Proceedings of the National Academy of Sciences of
  the United States of America}, vol. 109, no.~27, pp. 11\,031--11\,036, 2012.

\bibitem{Nakahira2015}
Y.~Nakahira, N.~Matni, and J.~C. Doyle, ``{Hard limits on robust control over
  delayed and quantized communication channels with applications to
  sensorimotor control},'' in \emph{Proc. IEEE CDC}, 2015, pp. 7522--7529.

\bibitem{Nakahira2018}
Y.~Nakahira, Q.~Liu, N.~Bernat, T.~Sejnowski, and J.~C. Doyle, ``{Theoretical
  foundations for layered architectures and speed-accuracy tradeoffs in
  sensorimotor control},'' in \emph{Proc. IEEE ACC}, 2019, pp. 509--514.

\bibitem{Liu2019}
Q.~Liu, Y.~Nakahira, A.~Mohideen, A.~Dai, S.~Choi, A.~Pan, D.~M. Ho, and J.~C.
  Doyle, ``{Experimental and educational platforms for studying architecture
  and tradeoffs in human sensorimotor control},'' in \emph{Proceedings of the
  American Control Conference}, 2019, pp. 483--488.

\bibitem{Huang2011}
Y.~Huang and R.~P. Rao, ``{Predictive coding},'' \emph{Wiley Interdisciplinary
  Reviews: Cognitive Science}, vol.~2, no.~5, pp. 580--593, 2011.

\bibitem{Xue2020}
\BIBentryALTinterwordspacing
A.~Xue and N.~Matni, ``{Data-Driven System Level Synthesis},'' 2020. [Online].
  Available: \url{http://arxiv.org/abs/2011.10674}
\BIBentrySTDinterwordspacing

\bibitem{Angelucci2003}
A.~Angelucci and J.~Bullier, ``{Reaching beyond the classical receptive field
  of V1 neurons: Horizontal or feedback axons?}'' \emph{Journal of Physiology
  Paris}, vol.~97, no. 2-3, pp. 141--154, 2003.

\bibitem{Perge2009}
J.~A. Perge, K.~Koch, R.~Miller, P.~Sterling, and V.~Balasubramanian, ``{How
  the Optic Nerve Allocates Space, Energy Capacity, and Information},''
  \emph{Journal of Neuroscience}, vol.~29, no.~24, pp. 7917--7928, 2009.

\end{thebibliography}
	\bibliographystyle{IEEEtran}
\end{document}